\documentclass[dvipsnames,x11names,conference]{ieeeconf}

\IEEEoverridecommandlockouts                              
\usepackage{algorithm}
\usepackage{algorithmic}
\usepackage{xcolor}
\usepackage{lipsum}
\usepackage{hyperref}
\usepackage{wrapfig}

\usepackage{graphics} 
\usepackage{epsfig} 
\usepackage{mathptmx} 
\usepackage{times} 
\usepackage{amsmath} 
\usepackage{amssymb}  
\usepackage{extpfeil}

\usepackage{tikz}
\usetikzlibrary{plotmarks}
\usepackage{silence}
\usepackage[compatibility=false]{caption}
\usepackage{subcaption}
\usepackage{pgfplots}
\usepgfplotslibrary{fillbetween}
\pgfplotsset{compat=1.18}
\usepackage{comment}
\usepackage{xcolor}
\usepackage{pifont}
\usepackage{pgfplots}
\usepackage{ulem}
\usepackage{color}
\usepackage{float}
\usepackage{stfloats}
\DeclareMathAlphabet{\mathcal}{OMS}{cmsy}{m}{n}

    \pgfplotsset{
        compat=1.16,
        layers/MyLayers/.define layer set={
            waybackgroundlayer, boundingboxlayer,
            axis background,
            pre main,axis grid,axis ticks,axis lines,axis tick labels,
            axis descriptions,axis,main,foreground,
        }{/pgfplots/layers/axis on top},
        set layers=MyLayers,
    }

\usepackage{xcolor}

\newtheorem{proposition}{\textbf{Proposition}}

\def\KL{D_\text{KL}}

\usepackage{graphics} 
\usepackage{epsfig} 
\usepackage{mathptmx} 
\usepackage{times} 
\usepackage{amsmath} 
\usepackage{amssymb}  
\usepackage{mathrsfs}

\usepackage{cite}

\newcommand{\field}[1]{\mathbb{#1}}
\def\stateOT{\mathcal{X}}
\def\clB{{\cal B}}
\def\clP{{\cal P}}
\def\clJ{{\cal J}}
\def\transpose{{\hbox{\tiny \textsf{T}}}}
\def\Re{\field{R}}
\def\eqdef{\mathbin{:=}}
\def\Expect{{\sf E}}
\def\Cov{\text{\rm Cov}\,}

\title{\LARGE \bf
Moment Constrained Optimal Transport for \\ Thermostatically Controlled Loads
}

\author{ \parbox{2 in}{\centering Thomas Le Corre\\
         Inria, Paris, France \\ DI ENS, ENS, PSL University \\ Paris, France\\
         {\tt\small thomas.le-corre@inria.fr}}
         \hspace*{ 0.3 in}
         \parbox{2 in}{ \centering Julien Cardinal\\
        Inria, Paris, France \\ DI ENS, ENS, PSL University \\ Paris, France\\
         {\tt\small julien.cardinal@inria.fr}}
         \hspace*{ 0.3 in}
         \parbox{2 in}{ \centering Ana Bušić\\
        Inria, Paris, France \\ DI ENS, ENS, PSL University \\ Paris, France\\
         {\tt\small ana.busic@inria.fr}}
}

\begin{document}

\maketitle
\renewcommand{\thefootnote}{} 
\footnotemark\footnotetext{This work was carried out in the framework of the AI-NRGY project, funded by France 2030 (Grant No: ANR-22-PETA-0004).}
\renewcommand{\thefootnote}{\arabic{footnote}} 
\thispagestyle{empty}
\pagestyle{empty}

\begin{abstract}
Controlling large populations of thermostatically controlled loads (TCLs), such as water heaters, poses significant challenges due to the need to balance global constraints (e.g., grid stability) with individual requirements (e.g., physical limits and quality of service). In this work, we introduce a novel framework based on Moment Constrained Optimal Transport (MCOT) for distributed control of TCLs. By formulating the control problem as an optimal transport problem with moment constraints, our approach integrates global consumption constraints and physical feasibility conditions into the control design.
This problem with high (or infinite) dimensionality can be reduced to a much lower finite-dimensional problem. The structure of this problem allows for computing the gradient with Monte Carlo methods by generating trajectories of TCLs. 
Contrary to all previous work, in our MCOT framework, it is possible to choose the sampling law, which considerably speeds up the calculations. 
This algorithm mitigates the need for extensive state-space discretization and significantly reduces computational complexity compared to existing methods. Numerical experiments in a water heater case study demonstrate that our MCOT-based method effectively coordinates TCLs under various constraints. We further extend our approach to an online setting, illustrating its practical applicability on simulated data. 

\end{abstract}

\section{INTRODUCTION}
The control of distributed systems such as \textit{Thermostatically Controlled Loads} (TCLs), requires methods capable of coordinating a large number of agents while simultaneously satisfying both global constraints (e.g. stability of the grid) and individual requirements (e.g. physical constraints, quality of service). The literature on demand-side management and distributed control has explored various strategies to balance these requirements. Techniques such as \textit{Demand Dispatch} \cite{hochberg2006demand} and ensemble control methods \cite{garabe2022ensemble} have been developed to address the challenges posed by fluctuations in energy demand. Furthermore, several recent studies have applied these techniques in diverse settings, ranging from electric vehicle (EV) charging control \cite{he2012optimal, sadeghianpourhamami2018ev} to power grid regulation \cite{amaraouali2021real, rezvanizaniani2014demand}.

Various stochastic optimal control methods  
have proven effective in managing large-scale systems while incorporating entropic regularization to ensure robustness \cite{klq, todorov2007linearly, chertkov2018stochastic}. These techniques are based on a probabilistic formulation, that facilitates the handling of a large number of agents, by the use of \textit{Mean Field Control} (MFC) approaches. In this framework, the number of agents is considered very large, so they can be modeled by the limit of their empirical distribution. The resulting optimization problem is generally simpler to solve. For TCLs, these types of mean-field models have long been considered \cite{laurent1994,malhame1985,malhame1988}.

\textit{Optimal Transport} (OT) has emerged as a fundamental tool for optimizing the cost of moving probability distributions \cite{villani2008optimal, peyre2019computational}. By introducing constraints on the moments of distributions, \textit{Moment Constrained Optimal Transport} (MCOT) \cite{alfonsi2021approximation} offers a novel perspective for addressing control problems in uncertain environments \cite{mcot}. This approach allows the implementation of various global constraints such as maximum total power consumption, consumption signal tracking, constraints on global consumption ramps,... At the same time, hard constraints of TCLs dynamics (temperature evolving according to an ODE, bounded temperature,...) are encoded directly in the problem, which ensures compliance with physical and quality of service constraints. The disadvantage of such methods is the size of the problem, which grows with the complexity of the agent and environment modeling, and can lead to difficulties in optimizing over long time horizons.

\textbf{Contributions: }The main contributions of this article are:
\begin{itemize}
    \item We formulate the problem of TCL control with large (and even infinite) state space as an MCOT problem. Discretizing the state space is no longer required, and Monte Carlo methods are used instead, 
     significantly reducing algorithmic complexity. A contribution of this article is that the structure of the problem allows for choosing the sampling law for the Monte Carlo method, which is particularly useful for speeding up the calculations.
    \item Our new MCOT framework allows one to use different sets of constraints for various control objectives, including signal tracking, consumption limit constraints, or ramp constraints. 
    These control objectives can be managed directly (i.e. without the need to design a specific tracking signal, which is needed in \cite{klq}).
    \item We extend this formulation to consider Model Predictive Control in the case study where the dynamic of (heterogeneous) water heaters (WHs) is stochastic, and the noise due to water usage is progressively discovered. 
    \item Finally, for the tracking objective, we extend the approach to the on-line signal tracking setting (i.e. when the tracking signal is revealed in real time). 
\end{itemize}
The paper is structured as follows. Section \ref{s:Prob} presents the problem of controlling WHs for a finite number of agents. Section \ref{s:MCOT} introduces the problem formulation as an MCOT problem, the main assumptions and the resolution of the dual problem. Section \ref{s:Alg} presents the algorithm and develops the Monte Carlo methods, and Section \ref{s:Num} illustrates the performance of our method in a concrete application. Finally, we introduce a Model Predictive Control method with an application on the dataset SMACH (\textit{Multi-agent Simulation of Human Activity in the Household}) \cite{SMACH} in Section \ref{s:On}.

\section{TCL CONTROL PROBLEM}
\label{s:Prob}
For the sake of clarity, this article will focus on the case of water heaters, but the application of this article is similar for other TCLs (fridges, air conditioners, etc.).

\subsection{Water Heaters}

We consider a large population of $N$ homogeneous \textit{Water Heaters} (WHs). At time t, the $i$-th WH is modeled by its mean temperature $\theta^i(t)\in \Theta$ and its power mode $m^i(t)\in \{0,1\}$ (Off/On). These WHs must then follow the \textit{Ordinary Differential Equation} (ODE):
\begin{equation}
\frac{d\theta(t)}{dt} = -\underbrace{\rho(\theta(t)-\theta_{amb})}_{\text{heat loss}} + \underbrace{\sigma m(t)p }_{\text{Joule effect}} - \underbrace{\sigma \epsilon(t)}_{\text{water drain}},
\label{e:ODE}
\end{equation}
with $\rho$ the fraction of heat loss by minute, $\sigma$ the specific heat capacity of the volume of water, $p$ the heating power, $\theta_{amb}$ the room temperature, and $\epsilon(t)$ the power equivalent of the water drains at time $t$.

We are interested in controlling their behavior over a day, discretized in $T=144$  intervals of size $\delta t = 10$ minutes. The previous ODE (\ref{e:ODE}) is thus discretized  with the following update:
\begin{equation}
\label{e:EqDiff}
\theta_{t+1} = \theta_t - (\rho (\theta_t-\theta_{amb}) + \sigma m_t p_{max} - \sigma \epsilon_t)\delta t
\end{equation}
where for $t\in\{0,\dots,T\}$, $\theta_t$, $m_t$ and $\epsilon_t$ are the values of $\theta$, $m$ and $\epsilon$ at time $t\times \delta t$. We will call the trajectory of the $i$th water heater the following:

$$X^i=\{(\theta_t,m_t)\}_{1\leq t\leq T}\in\stateOT=(\Theta\times\{0,1\})^{T}.$$

\begin{figure}[ht]
\vspace{0.1cm}
    \centering
    \begin{tikzpicture}[scale=0.45]
    \begin{axis}[xtick  = {0,4,8,12,16,20,24},legend style={at={(3,2.05)}, anchor=north east},legend columns=1,
xlabel={Time of day (h)},x label style ={at={(1.1,-0.33)},anchor=north},
grid=major,
width=0.7\textwidth,height=0.6\textwidth,
ylabel={Temperature (°c)},y label style ={at={(-0.55,0.25)},anchor=north west},
legend entries={$\theta$,$\theta_{min}$,$\theta_{max}$}],
    \addplot [draw=blue,ultra thick] table[x index=0,y index=1]{DataFigures/one_traj.txt};
    \addplot [draw=green,ultra thick] table[x index=0,y index=1]{DataFigures/Tmin.txt};
    \addplot [draw=red,ultra thick] table[x index=0,y index=1]{DataFigures/Tmax.txt};
    \end{axis}
    \end{tikzpicture}
    \caption{Example of a nominal trajectory of a single water heater, starting at $\theta_0=54$°c and $m_0=0$}
    \label{fig:NomTraj}
    \vspace{-0.3cm}
\end{figure}
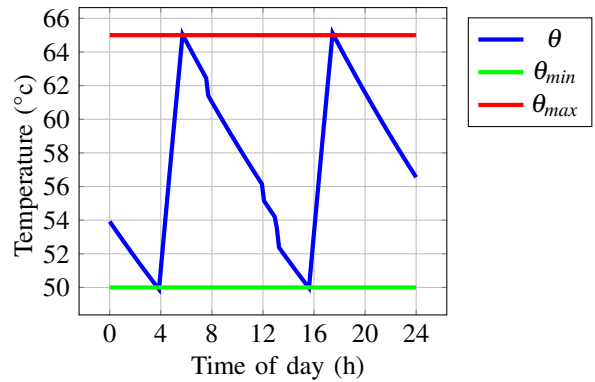

We define the nominal dynamics as the uncontrolled dynamic \cite{bianca,klq}. A water heater aims at keeping its mean temperature between $\theta_{min}=50^{\circ}C$ and $\theta_{max}=65^{\circ}C$ by turning the water heater Off whenever the temperature reaches  $\theta_{max}$ and turning it back On whenever the temperature reaches below $\theta_{min}$. A nominal trajectory is represented in Fig. \ref{fig:NomTraj}.

\subsection{Problem for N Water Heaters}

We consider a central agent who aims to coordinate $N$ WHs to satisfy $A$ global constraints on the aggregate behavior of the WHs:
\begin{equation}
\label{e:ConstraintsN}
    \forall a\in\{1,\dots,A\}, \ \sum\limits_{i=1}^N f^{(a)}(X^i)\leq0,
\end{equation}
where the function $f$ is defined as $f\colon\stateOT\to\Re^A$.

An example of such a global constraint would be a global max power $P_{max}$ to be respected at every moment of the day. In this case, there are $A=T$ constraints: $\forall a, f^{(a)}(X^i)=m_ap-\frac{P_{max}}{N}$.
Equality constraints can be expressed as pairs of inequality constraints, so it is possible to impose equality constraints when needed, such as tracking a signal for the total consumption for example. These and other examples of constraints are illustrated in Section \ref{s:Num}.
To achieve its objective of respecting the constraints, the central agent can change the trajectories of the water heaters from $\{X^i\}_{1\leq i\leq N}$ to $\{Y^i\}_{1\leq i\leq N}$, but they must respect the ODE defined above and have the same state at $t=0$. However, this change induces some cost $c$. This leads to the following problem:

\begin{equation}
\label{e:ProbN}
\begin{aligned}
    \min\limits_{\{Y^i\}_{1\leq i\leq N}} \Bigl\{ &\sum\limits_{i=1}^N c(X^i,Y^i) 
    : \sum\limits_{i=1}^K f(Y^i)\leq 0, \\
    &\forall i, X^i_0=Y^i_0 \text{ and } \forall i, Y^i \text{ follows ODE } \eqref{e:ODE} \Bigr\}
\end{aligned}
\end{equation}

where $X^i$ are the nominal trajectories.

\section{MOMENT CONSTRAINED OPTIMAL TRANSPORT}
\label{s:MCOT}
We formulate Problem (\ref{e:ProbN}), as a Mean Field Control problem, using the \textit{Moment Constrained Optimal Transport} (MCOT) approach. The dual of this formulation is computed, and an explicit expression of its gradient is obtained.

\subsection{Mean Field Formulation}

In this article, we focus on the Mean Field formulation \cite{lasry2007mean} of the Problem (\ref{e:ProbN}) (which intuitively corresponds to $N\rightarrow\infty$). The set of agents $\{X^i\}_{1\leq i\leq N}$ is thus considered through their common probability distribution $\mu\in\clB(\stateOT)$, with $\clB(\stateOT)$ the set of Borel probability measures on $\stateOT$.
This formulation is based on the assumption that the agents are homogeneous. In Section \ref{s:Num}, we will see how this assumption can be relaxed to accommodate heterogeneous water heaters within this framework.
We denote the nominal distribution $\mu_1$ corresponding to the nominal dynamics defined previously and illustrated in Fig. \ref{fig:NomTraj}.

Changing the trajectory from $X^i$ to $Y^i$ in Problem (\ref{e:ProbN}) becomes a transport of mass between two distributions of trajectories in this Mean Field formulation. We will therefore use the tools of optimal transport. Our objective is to find a bivariate distribution $\pi\in\mathcal{B}(\stateOT\times\stateOT)$, where $\mathcal{B}(\stateOT\times\stateOT)$ is the set of Borel probability measures on $\stateOT\times\stateOT$. 

\subsection{Problem Formulation}

A bivariate distribution $\pi$ is associated to a cost:

$$ \langle \pi,c\rangle=\int c(x,y)\pi(x,y)dxdy.$$

We consider the marginals, defined by $\pi_1=\int\pi(x,y)dy$ and $\pi_2=\int\pi(x,y)dx$. The first marginal is the nominal distribution: $\pi_1=\mu_1$. The second marginal $\pi_2$ must respect the global constraints introduced in (\ref{e:ConstraintsN}) and thus belongs to the set:
\begin{equation}
\clP_f =\{   \mu\in \clB(\stateOT)  : \forall   a\in\{1,\dots,A\},  \langle \mu,f^{(a)}\rangle \leq 0  \     \} .
\label{e:SimplexConstrained}
\end{equation}

As in Problem (\ref{e:ProbN}), we do not want to modify the distribution at time $0$, as it is considered part of the problem, and we try to control the rest of the trajectory. Thus, the  bivariate distribution $\pi$ must belong to
\begin{align*}
    \mathcal{U}(\mu_1)=\{& \pi \in \mathcal{B}(\stateOT\times\stateOT) : \pi_1=\mu_1 \text{ and } \forall x,y\in\stateOT, \\&\pi((x_0,x_{1:T}),(y_0,y_{1:T}))=0  \text{ if } x_0\neq y_0 \},
\end{align*}
where we denote a trajectory $x= (x_0,x_{1:T})\in\stateOT$ with $x_0$ its state at time $0$ and $x_{1:T}$, the rest of the trajectory.

We introduce the following entropic regularizer as in \cite{cuturi2013sinkhorn}.
\begin{equation}
\KL (\pi \| \mu_1\otimes\mu_2 ) = \int_{\stateOT\times\stateOT} \log\Bigg(\frac{\pi(x,y)}{\mu_1(x)\mu_2(y)}\Bigg) \pi(x,y)dxdy,
\label{e:OurReg}
\end{equation}
The role of this regularizer is twofold: (i) for computational reasons, to obtain an explicit expression of the gradient of the dual problem (ii) This term is infinite if the support of $\pi_2$ is larger than the support of $\mu_2$. So by choosing $\mu_2$ null for trajectories that are not following the ODE, we can ensure that this ODE will be verified for $\pi_2$.

We can introduce the Mean Field Control problem, studied in this paper, as a Mean Field formulation of Problem (\ref{e:ProbN}):

\textbf{Problem MCOT-C: \it Moment Constrained Optimal Transport for Control}
\begin{equation}
\label{e:MCOTC}
    \min\limits_{\pi}  \bigl\{ \langle \pi , c \rangle + \varepsilon \KL (\pi\|\mu_1\otimes\mu_2) :   \pi \in  \mathcal{U}(\mu_1)\,,  \  \pi_2 \in \clP_f \bigl\}.
\end{equation}
As the $\KL$ term is used as a regularizer (its value is not what we want to minimize), low values of $\varepsilon$ will be considered.

\subsection{Dual problem}

\textbf{Assumptions}

The following assumptions are introduced for the existence of optimizers and desirable properties of the dual:

\textbf{(A1)}
$c\colon \mathcal{X}\times \mathcal{X} \to\Re_+$ and  $f\colon\mathcal{X}\to\Re^A$ are continuous, and there is an open neighborhood $\mathcal{N}\subset \mathbf{R}^A$ containing $0$ such that
$\mathcal{P}_{f,r}=\{   \mu\in \clB(\stateOT)  : \forall   a\in\{1,\dots,A\},  \langle \mu,f^{(a)}\rangle \leq r_a  \     \} $ is non-empty for all $r\in \mathcal{N}$.

\textbf{(A2)} 
$\mu_1$ and $\mu_2$ have compact support, and the problem is feasible under perturbations:   for any $r\in \mathcal{N}$,  
 there is $\pi$ satisfying  $ \pi_2 \in \mathcal{P}_{f,r}$ and 
 $\pi \in  \mathcal{U}(\mu_1)$.
   
\textbf{(A3)} 
 $\Sigma^0 \eqdef \Cov(Y)$ is positive definite when $Y\sim\mu_2$.
 
\textbf{Dual}
The dual of MCOT-C is by definition
the function 
\begin{equation}
\label{e:DualMCOTC}
\begin{aligned}
    \varphi^*&(\lambda)=
\min\limits_{\pi}  \bigl\{ \varepsilon\KL (\pi\|\mu_1\otimes\mu_2) -\langle \pi , \ell_0^\lambda \rangle  :   \pi \in \mathcal{U}(\mu_1)\},
\end{aligned}
\end{equation}
where $\lambda$ is the Lagrange multiplier associated with the constraints in $\mathcal{P}_f$ and $\forall x,y\in  \stateOT,\ell_0^\lambda(x,y)  =  -\lambda^\transpose f(y) - c(x,y)  \,$.

For each $\lambda\in\Re_+^A$, $\varepsilon>0$ and $x=(x_0,x_{1:T})\in  \stateOT$, we denote 
\begin{equation}
B_{\lambda,\varepsilon}(x) = \varepsilon \log \int \exp\bigl(   \varepsilon^{-1}\ell_0^\lambda((x_0,x_{1:T}),(x_0,y_{1:T})  \bigr)\mu_2(y)dy_{1:T}   .
\label{e:RegMGF}
\end{equation}

\begin{subequations}

\begin{proposition}
\label{t:MCOT-C}
\textbf{(i)}
The minimum \eqref{e:DualMCOTC} 
gives:
$$
\varphi^*( \lambda)  =   -   \langle \mu_1,  B_{\lambda,\varepsilon} \rangle.
$$

\textbf{(ii)}
The maximizer is $\pi^\lambda(x,y)=P^\lambda(x,y)\mu_1(x)$ with  $\forall x=(x_0,y_{1:T})\in\stateOT, \forall y=(x_0,y_{1:T})\in \stateOT$
\begin{equation}
 \label{e:TlambdaFPR}
    P^\lambda(x,y) = \mu_2(y) \delta_{x_0}(y_0)  \exp(L^\lambda(x,y))  \,, 
\end{equation}
\begin{equation}      
 L^\lambda(x,y)   = \varepsilon^{-1}  \{ 	l_0^\lambda(x,y) -  B_{\lambda,\varepsilon}(x) \} \,, 
\end{equation}

\textbf{(iii)}    
There is no duality gap:  there is a unique $\lambda^*\in\Re_+^A$ satisfying
\begin{equation}
\begin{aligned}
\varphi^*&( \lambda^*)  =\\ &\min\limits_{\pi}  \bigl\{ \langle \pi , c \rangle + \varepsilon \KL (\pi\|\mu_1\otimes\mu_2) :   \pi \in  K(\mu_1)\,,  \  \pi_2 \in \clP_f \bigl\},
\end{aligned}
\label{e:NoGap1S-RMCOT}
\end{equation} 
\end{proposition}
\end{subequations}
\begin{proof} The proof is based on convex duality between relative entropy and log moment generating functions. For any probability measure $\mu$ on $\stateOT$  and function   $g\colon\stateOT\to\Re$,  the log moment generating function is denoted $\Lambda_{\mu}(g) = \log \langle \mu ,  e^g \rangle$. The following result is a standard tool in information theory \cite{demzei98a}.
\begin{subequations}
 For Borel measurable   $g\colon\stateOT\to\Re$,
\begin{equation*}
 \Lambda_\mu (g) =  \sup_p  \{  \langle p , g\rangle - \KL (p\| \mu)  \}.
\label{e:Lambda2D}
\end{equation*}
If  $\Lambda_\mu(g)<\infty$ then the supremum is achieved, where the optimizer $p^*$ has log likelihood ratio  $ \log(dp^*/d\mu)  = g - \Lambda_\mu(g)$.
\end{subequations}
By writing with an abuse of notation that
\begin{equation}
\begin{aligned}
    \varphi^*&(\lambda)=\\&-\varepsilon\langle\mu_1,\sup\limits_{P(x,.)}  \bigl\{  \langle P(x,.) , \varepsilon^{-1}\ell_0^\lambda(x,.) \rangle -\KL (P(x,.)\|\mu_2) \}\rangle,
\end{aligned}
\end{equation}
we can use this result and obtain the proposition.
\end{proof}

We make the convenient change of variables $\zeta = \varepsilon^{-1} \lambda$,   and consider $$ \clJ(\zeta)  \eqdef - \varepsilon^{-1}  \varphi^*(\varepsilon\zeta ). $$ 

We turn next to the representation of the derivatives of the dual function.   
The quantity $ \varepsilon^{-1}   B_{\varepsilon\zeta,\varepsilon}(x)$  is a log moment generating function for each $x$;  for this reason, it is not difficult to obtain suggestive expressions for the first and second derivatives with respect to $\zeta$.
\begin{subequations} 
\begin{proposition}
\label{t:RegDualCalculus} 
The function $\clJ$ is convex and continuously differentiable.   
The first and second derivatives of $\clJ$ admit the following representations:
\begin{equation} 
\nabla  \clJ(\zeta)   =    \langle \mu^\lambda , f\rangle   \,, 
\end{equation}
\vspace{-0.25cm}
\begin{equation}
\nabla^2  \clJ(\zeta)  = \Expect^\lambda[f(Y)  f(Y) ^\transpose - \Expect^\lambda[ f(Y) \mid X] \Expect^\lambda[f(Y) \mid X] ^\transpose   \bigr],
\label{e:GradHessRegCov}
\end{equation} 
where $\mu^\lambda=\pi_2^\lambda$. It follows that $\clJ $ is strictly convex.
\end{proposition}
\end{subequations}
\begin{proof}
    The proof is straightforward. We simply derive the equations obtained in Proposition \ref{t:MCOT-C}.
\end{proof}

\section{ALGORITHM}
\label{s:Alg}

\subsection{Monte Carlo Methods}
Recall the expression for $\mu^\lambda$, defined as the second marginal of $\pi^\lambda$ obtained in proposition \ref{t:MCOT-C}, (\ref{e:TlambdaFPR}): $\forall y=(y_0,y_{1:T})\in \stateOT$,

\begin{equation}
    \label{e:muL}
    \mu^\lambda(y)=\mu_2(y) \int_x\delta_{x_0}(y_0)\exp(L^\lambda(x,y))\mu_1(x)dx
\end{equation}

$\mu^\lambda(y)$ is thus proportional to $\mu_2(y)$, which leads to two key points of this article:

\begin{enumerate}
    \item The support of $\mu^\lambda$ is included in the support of $\mu_2$. We can thus impose zero values on $\mu^\lambda$ by choosing $\mu_2$ to be zero on these values. This is particularly useful in control applications where certain states (or trajectories) are not physically possible. In the water heater example, we use this point to have only trajectories respecting the differential equation (\ref{e:EqDiff}). Any other trajectory is linked to a zero value for $\mu_2$ and therefore for $\mu^\lambda$.
    \item The term $w^\lambda(y)=\int_x\delta_{x_0}(y_0)\exp(L^\lambda(x,y))\mu_1(x)dx$ can be interpreted as a weight. Thus, $\mu^\lambda$ is a re-weighting of $\mu_2$, which leads to Importance Sampling-type methods \cite{wasserman2010}. An appropriate choice of $\mu_2$ can therefore reduce the complexity of the gradient computation at each stage.
\end{enumerate}

We can generate $Z$ trajectories $\{Y_z\}_{1\leq z\leq Z}$ according to $\mu_2$ and compute the gradient, with the following approximation:

\begin{equation}
\label{e:MCApprox}
    \langle \mu^\lambda , f\rangle\approx\frac{1}{Z}\sum\limits_{z=1}^Z w^\lambda(Y_z)f(Y_z).
\end{equation}

\begin{proposition}
\label{t:IS}
(i) The approximation in (\ref{e:MCApprox}) has zeros bias and:  $$\frac{1}{Z}\sum\limits_{z=1}^Z w^\lambda(Y_z)f(Y_z) \xlongrightarrow[\ Z \to \infty\ ]{\text{a.s.}} \langle \mu^\lambda , f\rangle,$$

(ii) The variance of this estimator is:
\begin{equation}
\label{e:Var}
    \frac{1}{Z}\int f^2(y)w^\lambda(y)\mu^\lambda(y)dy-\frac{1}{Z}\Big(\int f(y)\mu^\lambda(y)dy\Big)^2,
\end{equation}
and to minimize this expression, we should have:
\begin{equation}
\label{e:MinG}
    \frac{\mu_2(y)}{\mu^\lambda(y)}=\frac{f(y)}{\int f(s)\mu^\lambda(s)ds}.
\end{equation}
\end{proposition}

\begin{proof}
    (i) is a straightforward computation and the convergence is given by the law of large numbers. For (ii), obtaining (\ref{e:Var}) is also straightforward and (\ref{e:MinG}) is obtained by proving that for $\frac{1}{w^\lambda(y)}=\frac{f(y)}{\int f(s)\mu^\lambda(s)ds}$, the variance is zero, thus minimal. 
\end{proof}

The Proposition \ref{t:IS} gives us an idea of how to choose $\mu_2$. First, (\ref{e:Var}) implies that we should choose $\mu_2$ to have thicker tails than $\mu^\lambda$ (if not, the term $w^\lambda(y)=\mu^\lambda(y)/\mu_2(y)$ may grow fast with $y$ and we may have a large and even infinite integral). For control, these integrals are mostly performed over a bounded set (e.g. set of temperatures for WHs), so this is not an issue. Equation (\ref{e:MinG}) implies that $\mu_2$ and  $\mu^\lambda$ should be close. In practice, when designing $\mu_2$, $\mu^\lambda$ is not known (it is the optimum we are looking for), but having a prior idea of the shape of $\mu^\lambda$, a good choice of $\mu_2$ might decrease the variance and thus improve the convergence of algorithm \ref{a:MCOT}.

Choosing the law from which trajectories are drawn to compute the gradient was not previously present in the literature on Mean Field Control for TCLs. For Kullback-Leibler Quadratic Optimal Control \cite{klq} in particular, which gradient resembles that of this article, calculating the gradient using Monte Carlo methods requires trajectories to be drawn according to the $\mu_1$ law. This raised 2 concerns: (i) some potentially desirable trajectories may have zero or very low probability in the $\mu_1$ distribution, and will never or almost never be drawn (ii) if the optimal distribution is far from $\mu_1$, the variance is large and it may take a long time to calculate a correct approximation of the gradient.

\subsection{Gradient Descent Algorithm}

We design the following Stochastic Gradient Descent (SGD) algorithm using the approximation of the gradient introduced in (\ref{e:MCApprox}). With $K$ the number of iterations, the step \(\{\rho_k\}_k\) is chosen empirically to accelerate the convergence.

\begin{algorithm}[ht]
\caption{SGD Algorithm for MCOT}
\label{a:MCOT}
\begin{algorithmic}[1]
\STATE Initialize \( \lambda_0 \)
\FOR{each \( k = 1,\dots, K \)}
    \STATE Generate $Z$ independent states realizations $(Y_1,...,Y_Z)$ according to $\mu_2$.
    \STATE $G_{k}\gets\frac{1}{Z}\sum_{z=1}^Z w^\lambda(Y_z)f(Y_z)$
    \STATE $ \lambda_{k+1}\gets\lambda_k-\rho_kG_k$
\ENDFOR
\end{algorithmic}
\end{algorithm}

\section{NUMERICAL RESULTS}
\label{s:Num}
\subsection{Presentation of the use case}

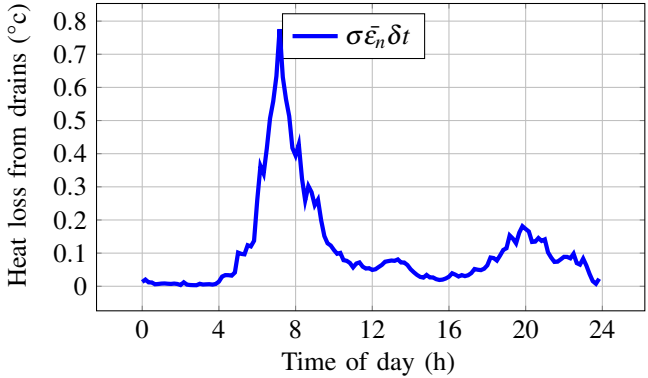
\begin{figure}[ht]
\vspace{0.1cm}
    \centering
    \begin{tikzpicture}[scale=0.45]
    \begin{axis}[xtick  = {0,4,8,12,16,20,24},legend style={at={(1.23,2.05)}, anchor=north east},legend columns=2,
xlabel={Time of day (h)},x label style ={at={(1.0,-0.33)},anchor=north},
grid=major,
width=1\textwidth,height=0.6\textwidth,
ylabel={Heat loss from drains (°c)},y label style ={at={(-0.5,0.0)},anchor=north west},legend entries={$\sigma\bar{\epsilon_n}\delta t$}],
    \addplot [draw=blue,ultra thick] table[x index=0,y index=1]{DataFigures/drains.txt};
    \end{axis}
    \end{tikzpicture}
    \caption{Average heat loss generated by the water drains during the day}
    \label{fig:DrainsDistrib}
\end{figure}

In this paper, we use simulated water drain data from \cite{bianca}, which was generated from the SMACH (\textit{Multi-agent Simulation of Human Activity in the Household}) platform \cite{SMACH} and averaged on 10-minute intervals. The distribution of the water drains $\epsilon_t$ throughout the day is represented in Fig. \ref{fig:DrainsDistrib}, with the average water drains noted $\bar{\epsilon_t}$. We separate the data between a training set and a validation set.

The full update equation for the nominal policy ($\mu_1$) can then be written as follows
\begin{equation}
\left \{
\begin{array}{l}
\theta_{t+1} = \theta_t - \rho\delta t (\theta_t-\theta_{amb}) + \sigma \delta t m_t p_{max} - \sigma \epsilon_t
\vspace{0.5em}\\
m_{t+1} = \left \{
\begin{array}{cc}
     m_t & \text{if } \theta_{t+1}\in[\theta_{min},\theta_{max}] \\
     0 & \text{if } \theta_{t+1} \geq \theta_{max} \\
     1 & \text{if } \theta_{t+1} \leq \theta_{min} \\
\end{array}
\right.
\vspace{0.5em}\\
\theta_0,m_0 \sim \nu_0
\end{array}
\right. 
\end{equation}
with $\nu_0$ being the initial density: the initial temperature $\theta_0$ is uniformly distributed in $[\theta_{min},\theta_{max}]$ and the initial mode $m_0$ is 1 (On) with probability $0.25$ and 0 (Off) with probability $0.75$

\begin{figure}[ht]
    \centering
    \begin{tikzpicture}[scale=0.45]
    \begin{axis}[xtick  = {0,4,8,12,16,20,24},legend style={at={(3,2.05)}, anchor=north east},legend columns=1,
xlabel={Time of day (h)},x label style ={at={(1.1,-0.33)},anchor=north},
grid=major,
width=0.7\textwidth,height=0.6\textwidth,
ylabel={Temperature (°c)},y label style ={at={(-0.55,0.25)},anchor=north west},
legend entries={$\theta$,$\theta_{min}$,$\theta_{max}$}],
    \addplot [draw=blue,ultra thick] table[x index=0,y index=1]{DataFigures/one_traj2.txt};
    \addplot [draw=green,ultra thick] table[x index=0,y index=1]{DataFigures/Tmin.txt};
    \addplot [draw=red,ultra thick] table[x index=0,y index=1]{DataFigures/Tmax.txt};
    \end{axis}
    \end{tikzpicture}
    \caption{Example of a controlled trajectory of a single water heater, starting at $\theta_0=54$°c and $m_0=0$ with two switches, one at $t_1=$1:00 and the other at $t_2=$14:20}
    \label{fig:2switchTraj}
\end{figure}
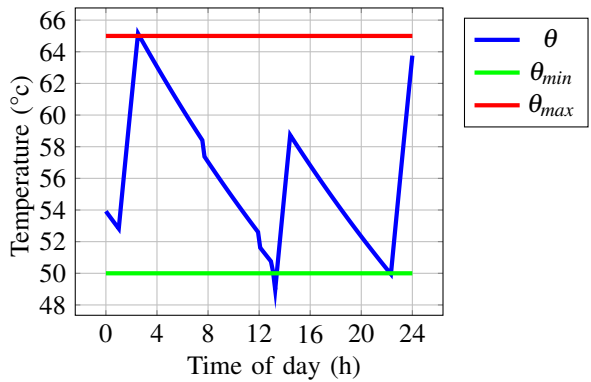

In the remainder of this paper, we allow the WHs to switch their power mode (from On to Off or from Off to On), while the temperature is still between the two bounds $\theta_{t+1}\in[\theta_{min},\theta_{max}]$. We will call switches those cases (i.e. the change of mode proposed by our model). We limit ourselves to two switches per day per water heater. We will note these two times $t_1,t_2\in\{1,\dots,T\}^2$. This limitation avoids frequent switching, which is undesirable for the water heater. Fig. \ref{fig:2switchTraj} shows an example of a trajectory with two switches.

Therefore, the power mode update for the policy $\mu_2$ can be written as follows:

\begin{equation}
m_{t+1} = \left \{
\begin{array}{cc}
     m_t & \text{if } \theta_{t+1}\in[\theta_{min},\theta_{max}] \text{ and } t\notin\{t_1,t_2\}\\
     (1-m_t) & \text{if } \theta_{t+1}\in[\theta_{min},\theta_{max}] \text{ and } t\in\{t_1,t_2\}\\
     0 & \text{if } \theta_{t+1} \geq \theta_{max} \\
     1 & \text{if } \theta_{t+1} \leq \theta_{min} \\
\end{array}
\right.
\end{equation}

\textbf{Tracking constraints:}
With $\{r_t\}_{1\leq t \leq T}\in \mathbb{R}^{T}$ the signal to track as the consumption of the whole population of WHs.
We can define the tracking constraint as $T \times 2$ inequality constraints:
\begin{equation}
\label{e:tracking}
\begin{array}{l}
\forall t\in\{1,\dots,T\},\\
\left \{
\begin{array}{l}
\sum\limits_{i=1}^N f_t(X^i)\leq0 \\
\sum\limits_{i=1}^N g_t(X^i)\leq0
\end{array}
\right.
\text{with }
\left \{
\begin{array}{l}
f_t(X^i) = (m_t^i-r_t/N)\\
g_t(X^i) = (r_t/N-m_t^i)\\
\forall i \in \{1,\dots,N\}
\end{array}
\right.
\end{array}
\end{equation}
\noindent
having $\forall i\in\{1,\dots,N\},\ X^i = \{(m_t^i,\theta_t^i)\}_{1\leq t \leq T} \in \mathcal{X}$

\textbf{Consumption limit constraints:}
With $\{u_t\}_{1\leq t \leq T}\in \mathbb{R}^{T}$ an upper bound on the consumption at each time step $t$.
We can define the upper bound constraint as $T$ inequality constraints:
\begin{equation}
\label{e:upperb}
\begin{array}{l}
\forall t\in\{1,\dots,T\},\\
\sum\limits_{i=1}^N f_t(X^i)\leq0 
\text{ with }
\left \{
\begin{array}{l}
f_t(X^i) = m_t^i-u_t/N\\
\forall i \in \{1,\dots,N\}
\end{array}
\right.
\end{array}
\end{equation}

\textbf{Ramp constraints:}
With $\{d^+_t\}_{1\leq t \leq T}\in \mathbb{R}^{T-1}$ an upper bound, and $\{d^-_t\}_{1\leq t \leq T}\in \mathbb{R}^{T-1}$ a lower bound on the gradient of the consumption at each time step $t$.
We can define the ramp constraint as $2\times (T-1)$ inequality constraints:
\begin{equation}
\begin{array}{l}
\forall t\in\{1,\dots,T-1\},\\
\left \{
\begin{array}{l}
\sum\limits_{i=1}^N f_t(X^i)\leq0 \\
\sum\limits_{i=1}^N g_t(X^i)\leq0
\end{array}
\right.
\text{with }
\left \{
\begin{array}{l}
f_t(X^i) = m_{t+1}^i-m_{t}^i - d^+_t/N\\
g_t(X^i) = m_{t}^i-m_{t+1}^i - d^-_t/N\\
\forall i \in \{1,\dots,N\}
\end{array}
\right.
\end{array}
\label{e:gradconstraints}
\end{equation}

We can combine several types of constraints and define other ones.\\
\textbf{Local Cost:}
Moreover, we can define local costs as a distance between a nominal trajectory and its controlled counterpart. 
In this paper, we focused on compliance with constraints, so we chose the simple cost function:
\begin{equation}
\begin{array}{cccc}
 c:&\mathcal{X}\times\mathcal{X} & \rightarrow  & \mathbb{R}\\
  & X,Y & \mapsto & 
  \left |
  \begin{array}{l}
    0 \text{ if } X = Y  \\
    1 \text{ otherwise }
  \end{array}
  \right.
\end{array}
\end{equation}

Depending on the use case, we can define alternative costs. For example, we can choose $c$ as the difference between the daily consumption of the nominal trajectory and the controlled one, to reduce the impact on user consumption. We can also choose $c$ as the number of switches added during the day or the temperature difference at the end of the day, to reduce the impact of the algorithm on user comfort.

\subsection{Results}
In this section, we evaluate the global consumption of the population of WHs as a percentage of the maximum global consumption. This way the global consumption $G$ at each time step is a value between $0$ and $1$, that can be expressed as follows:
\begin{equation}
\forall t\in \{1,\dots,T\}\ G_t = \frac{1}{N}\sum\limits^N_{i=1}m_t^i.
\end{equation}

We will refer to $G$ as the consumption of the population of WHs. To compute the results, unless stated otherwise, we will consider a population of $N=2000$ WHs.

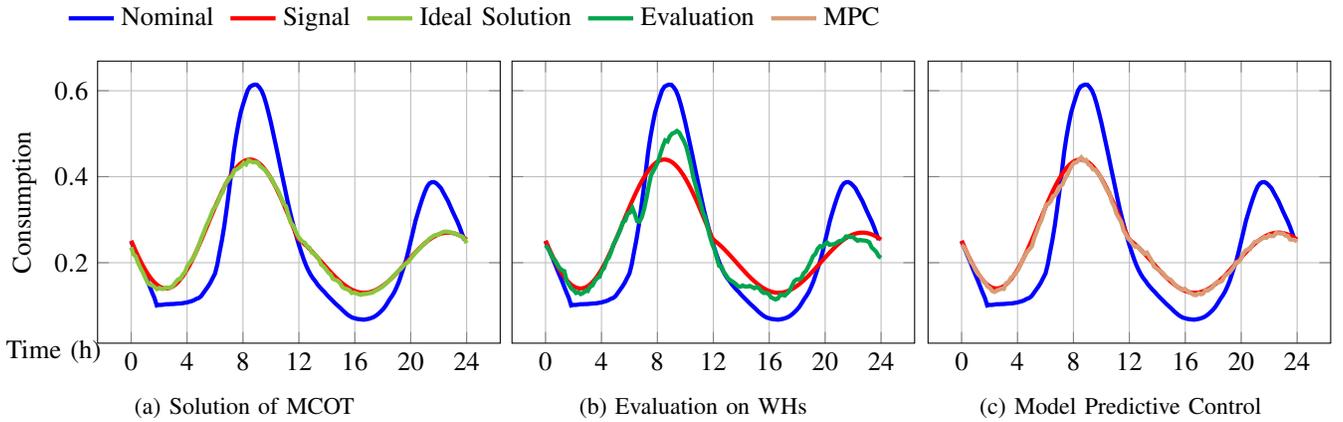
\begin{figure*}[t!]
\captionsetup[subfigure]{}
\hspace{-0.35cm}
\subcaptionbox{Solution of MCOT}{
\begin{tikzpicture}[scale=1]
\begin{axis}
[xtick  = {0,4,8,12,16,20,24},
xlabel={Time (h)},x label style ={at={(-0.25,0.045)},anchor=north west},
grid=major,ylabel={Consumption},
width=0.39\textwidth,height=0.3\textwidth]
\addplot [draw=blue,ultra thick] table[x index=0,y index=1]{DataFigures/nominal.txt};
\addplot [draw=red,ultra thick] table[x index=0,y index=1]{DataFigures/Signal_to_track.txt};
\addplot [draw=LimeGreen,ultra thick, name path=f] table[x index=0,y index=1]{DataFigures/Mean_Field_tracking.txt};
\path [name path=xaxis]
      (\pgfkeysvalueof{/pgfplots/xmin},0) --
      (\pgfkeysvalueof{/pgfplots/xmax},0);
\end{axis}
\end{tikzpicture}}%
\hspace{-6cm}
\subcaptionbox{Evaluation on WHs}{\begin{tikzpicture}[scale=1]
\begin{axis}
[yticklabel=\empty,legend style={at={(0.45,1.23)}, anchor=north,draw=none},legend columns=5,grid=major,xtick  = {0,4,8,12,16,20,24},
width=0.39\textwidth,height=0.3\textwidth,
legend entries={Nominal \quad , Signal \quad , Ideal Solution \quad,  Evaluation \quad, MPC \space\space\space\space\space\space\space\space\space\space\space\space\space\space\space\space\space\space\space\space\space\space\space\space\space\space\space\space\space\space\space\space\space\space\space\space\space\space\space\space\space\space\space\space\space\space\space}],
\addplot [draw=blue,ultra thick] table[x index=0,y index=1]{DataFigures/nominal.txt};
\addplot [draw=red,ultra thick] table[x index=0,y index=1]{DataFigures/Signal_to_track.txt};
\addlegendimage{LimeGreen,ultra thick}
\addplot [draw=Green,ultra thick, name path=f] table[x index=0,y index=1]{DataFigures/offline_first_step_learned.txt};
\addlegendimage{Tan,ultra thick}
\path [name path=xaxis]
      (\pgfkeysvalueof{/pgfplots/xmin},0) --
      (\pgfkeysvalueof{/pgfplots/xmax},0);
\end{axis}
\end{tikzpicture}}
\hspace{-5.8cm}
\subcaptionbox{Model Predictive Control}{\begin{tikzpicture}[scale=1]%
\begin{axis}
[grid=major,yticklabel=\empty,xtick  = {0,4,8,12,16,20,24},
width=0.39\textwidth,height=0.3\textwidth],
\addplot [draw=blue,ultra thick] table[x index=0,y index=1]{DataFigures/nominal.txt};
\addplot [draw=red,ultra thick] table[x index=0,y index=1]{DataFigures/Signal_to_track.txt};
\addplot [draw=Tan,ultra thick, name path=f] table[x index=0,y index=1]{DataFigures/offline_tracking.txt};
\path [name path=xaxis]
      (\pgfkeysvalueof{/pgfplots/xmin},0) --
      (\pgfkeysvalueof{/pgfplots/xmax},0);
\end{axis}
\end{tikzpicture}}
\vspace{-0.15cm}
\caption{Aggregated consumption of a population of 2000 WHs on the tracking problem; (a) The solution in the ideal case; (b) Its evaluation on a finite population of WHs; (c) The MPC version by generating (a) and (b) every $\delta t = 10$ minutes.}
\vspace{-0.2cm}
\label{fig:tracking}
\end{figure*}

Using the drains from the SMACH dataset, with the average values at each time step drawn in Fig. \ref{fig:DrainsDistrib}, the nominal consumption is the consumption of the population of WHs following nominal trajectories.
The nominal consumption is shown in blue in Fig. \ref{fig:tracking} and reflects the drain distribution. The peak in water drains around 7:00 is responsible for the need of most WHs to turn On directly afterwards, to maintain their temperature above $50^{\circ}C$. This synchronization, accompanied by an increase in water drains in the evening is responsible for the second consumption peak at the end of the day.

In the tracking problem, we use the constraints in (\ref{e:tracking}) with $\{r_t\}_{1\leq t \leq T}\in [0,1]^{T}$ set to attenuate consumption peaks and to have the same average consumption as the nominal signal.

To evaluate our approach for a population of WHs we use the algorithm \ref{a:Evaluation}.
The average consumption at the end $G$ gives us the evaluation curve in Fig. \ref{fig:tracking}(b).

\begin{algorithm}[ht]
\caption{Training and Evaluation of MCOT for WHs}
\label{a:Evaluation}
\begin{algorithmic}[1]
\STATE Parameters: Initial WHs states $\{X_{0}^i\}_{1\leq i \leq N}$
\STATE --- Training MCOT on WHs ---
\STATE Initialize \( \lambda_0 \)
\FOR{each \( k = 1,\dots, K \)}
    \FOR{each \( i = 1, \dots, N \)}
        \STATE Generate $Z$ trajectories $(Y_1^i,...,Y_Z^i)$ with \\
        $\forall z \in \{1,\dots,Z\},\ Y_{0,z}^i=X_{0}^i$ according to $\mu_2$.
        \STATE $G_{k}^i\gets\frac{1}{Z}\sum_{z=1}^Z w^\lambda(Y_z^i)f(Y_z^i)$
    \ENDFOR
    \STATE $G_{k} = \frac{1}{N} \sum_{i=1}^N G_{k}^i$
    \STATE $ \lambda_{k}\gets\lambda_{k-1}-\rho_kG_k$ \label{a:Evaluation:lbda}
\ENDFOR
\STATE --- Evaluating MCOT on WHs ---

\FOR{each \( i = 1, \dots, N \)}
    \STATE Generate $Z$ trajectories $(Y_1^i,...,Y_Z^i)$ with \\
    $\forall z \in \{1,\dots,Z\},\ Y_{0,z}^i=X_{0}^i$ according to $\mu_2$.
    \STATE Draw $Y^i$ with the weights $\{w^\lambda(Y_1^i),\dots,w^\lambda(Y_Z^i)\}$
    \STATE Generate $\tilde{Y}^i$ with the same 2 (or less) added switches as $Y^i$ on the validation water drains
\ENDFOR
\STATE $G = \frac{1}{N} \sum_{i=1}^N f(\tilde{Y}^i)$
\end{algorithmic}
\end{algorithm}

\subsection{Comparison with Kullback Leibler Quadratic Optimal Control for heterogeneous Water Heaters}

In the heterogeneous tracking problem, the characteristics of the WHs differ. This problem can be solved with an easy extension of MCOT by computing the temperature update of each WH based on its own characteristics, during training and evaluation.

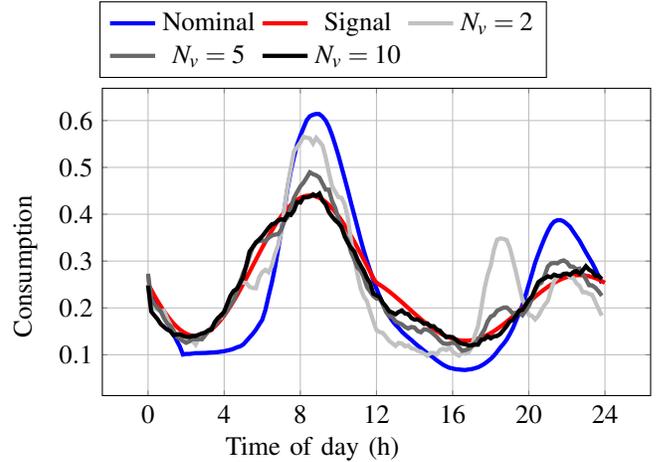
\begin{figure}[ht]
    \centering
    \begin{tikzpicture}[scale=0.45]
    \begin{axis}[xtick  = {0,4,8,12,16,20,24},legend style={at={(1.7,2.75)}, anchor=north east},legend columns=3,
xlabel={Time of day (h)},x label style ={at={(0.75,-0.33)},anchor=north},
grid=major,
width=1\textwidth,height=0.6\textwidth,
ylabel={Consumption},y label style ={at={(-0.5,0.25)},anchor=north west},legend entries={Nominal,Signal,$N_v = 2$,$N_v = 5$,$N_v = 10$}],
    \addplot [draw=blue,ultra thick] table[x index=0,y index=1]{DataFigures/nominal.txt};
    \addplot [draw=red,ultra thick] table[x index=0,y index=1]{DataFigures/Signal_to_track.txt};
    
    \addplot [draw=white!30!gray!70,ultra thick] table[x index=0,y index=1]{DataFigures/klq2.txt};
    \addplot [draw=gray!30!black!70,ultra thick] table[x index=0,y index=1]{DataFigures/klq5.txt};
    \addplot [draw=black,ultra thick] table[x index=0,y index=1]{DataFigures/klq10.txt};
    \end{axis}
    \end{tikzpicture}
    \caption{Comparison of the KLQ aggregated consumption, in evaluation, on average drains. We chose $N_v$ points of discretization of the volume range $[0.1,0.3]$ m$^3$.}
    \label{fig:KLQ}
    \vspace{-0.5cm}
\end{figure}

In this sub-section, we compare ourselves with an existing method in the literature \textit{Kullback Leibler Quadratic Optimal Control} (KLQ), on a signal tracking use case. The aim here is to demonstrate the interest of MCOT in cases where the space is very large, as there is no need to discretize the space to calculate the gradient, unlike methods such as KLQ. The use case chosen here is that of a population of water heaters with heterogeneous volumes. This is a difficult case to deal within the mean-field literature, and many methods like KLQ can handle it by grouping the heaters into “classes” and thus discretizing this heterogeneity \cite{KLQHeterogeneous}. There is no need to do this with MCOT, as trajectories can be generated with this heterogeneity. In order to compare these two models, we use the nominal behavior of KLQ, with a modified nominal distribution (the probability of switching increases as we approach the temperature limits) and we use KLQ in the context of perfect signal tracking (the quadratic term of KLQ is chosen with a very large parameter to strongly penalize the distance to the signal). For a lower tracking error (norm of the distance to the signal), KLQ requires strong discretization, as can be seen in Fig \ref{fig:KLQ}, but it takes more and more time and memory. The performance comparison shown in Table \ref{tab:1} shows that MCOT is faster and requires less memory to achieve a lower tracking error.

\begin{table}[ht]
    \centering
    \begin{tabular}{|c||c|c|c|c|}
        \hline
        Algorithm & $N_v$ & Execution Time & Memory & Error   \\
        \hline
        MCOT & - & 8.2 s & 1.66 GB& 0.111\\
        KLQ & 2 & 21.4 s & 4.96 GB& 0.823\\
        KLQ & 5 & 52.3 s & 9.75 GB& 0.294\\
        KLQ & 10 & 106.6 s & 17.75 GB& 0.2026\\
        \hline
    \end{tabular}
    \caption{Comparison between MCOT and KLQ for different discretizations. Algorithms run on a single GPU RTX6000 and evaluated on 2000 WHs with average drains. For KLQ, we discretize the temperature in 500 points.}
    \label{tab:1}
    \vspace{-0.5cm}
\end{table}

\section{MODEL PREDICTIVE CONTROL}
\label{s:On}
\subsection{The Model Predictive Control algorithm}

As we can see in Fig. \ref{fig:tracking}(b), the evaluation is good only in the short term and not good enough to precisely track the target signal throughout the day. This is why we introduce the Model Predictive Control (MPC) algorithm.

To better fit the signal, new decisions are made every $\delta t = 10$ minutes. These decisions take into account the current state of the population, and can then re-adjust if the stochastic water drains at the previous time step caused the consumption to overshoot or undershoot the signal to follow.

Algorithm \ref{a:MPC} does just that and keeps track of the remaining control switches for every WH so that no WH exceeds two switches inside the temperature limits.

\begin{algorithm}[ht]
\caption{Model Predictive Control for WHs}
\label{a:MPC}
\begin{algorithmic}[1]
\STATE Parameters: Initial WHs states $X_0$
\STATE Initialize $G^{MPC}$
\FOR{each \( t = 0, \dots, T \)}
    \STATE $G,\tilde{Y} \gets $ Algorithm \ref{a:Evaluation} initialized with $\{X_t^i\}_{1\leq i \leq N}$ \label{a:MPC:Evaluation}
    \STATE $X_{t+1} \gets \tilde{Y}_{t+1}$
    \STATE $G^{MPC}_t \gets G_t$
    \STATE Reduce the number of available switches for every water heater that switched at time step t
\ENDFOR
\end{algorithmic}
\end{algorithm}

\subsection{Results}

\begin{figure}[ht]
    \centering
    \begin{tikzpicture}[scale=0.45]
    \begin{axis}[xtick  = {0,4,8,12,16,20,24},legend style={at={(1.57,2.5)}, anchor=north east},legend columns=3,
xlabel={Time of day (h)},x label style ={at={(0.75,-0.33)},anchor=north},
grid=major,
width=1\textwidth,height=0.6\textwidth,
ylabel={Consumption},y label style ={at={(-0.5,0.25)},anchor=north west},legend entries={Nominal,Signal,MPC}],
    \addplot [draw=blue,ultra thick] table[x index=0,y index=1]{DataFigures/nominal_het.txt};
    \addplot [draw=red,ultra thick] table[x index=0,y index=1]{DataFigures/Signal_to_track.txt};
    \addplot [draw=Tan,ultra thick, name path=f] table[x index=0,y index=1]{DataFigures/offline_tracking_het.txt};
    
    \end{axis}
    \end{tikzpicture}
    \caption{Heterogeneous case, with the volume (m$^3$), height (m), isolation thickness (m), heat resistance power (W) of the water heaters, uniformly distributed in $[0.1,0.3] \times [0.8,2] \times [0.02,0.05] \times [1500,2900]$.}
    \label{fig:Heterogeneous}
\end{figure}
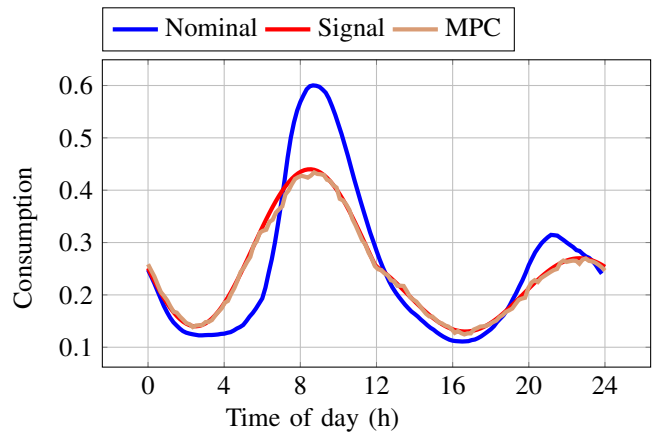

For the tracking problem, MPC works well in both the homogeneous case (Fig. \ref{fig:tracking}(c)) and the heterogeneous case (Fig. \ref{fig:Heterogeneous}(c)).

\textbf{In the homogeneous case}, we evaluate the MPC algorithm on multiple constrained settings.

\begin{figure*}[b]
\captionsetup[subfigure]{}
\hspace{-0.35cm}
\subcaptionbox{Whole day constraint}{
\begin{tikzpicture}[scale=1]
\begin{axis}
[xtick  = {0,4,8,12,16,20,24},
xlabel={Time (h)},x label style ={at={(-0.25,0.045)},anchor=north west},
grid=major,ylabel={Consumption},
width=0.39\textwidth,height=0.3\textwidth]
\addplot [draw=blue,ultra thick] table[x index=0,y index=1]{DataFigures/nominal.txt};
\addplot [draw=Tan,ultra thick, name path=f] table[x index=0,y index=1]{DataFigures/MPC_single_constraint.txt};
\addplot [draw=red,ultra thick, name path=f] table[x index=0,y index=1]{DataFigures/first_constraint.txt};
\path [name path=xaxis]
      (\pgfkeysvalueof{/pgfplots/xmin},0) --
      (\pgfkeysvalueof{/pgfplots/xmax},0);
\end{axis}
\end{tikzpicture}}%
\hspace{-6cm}
\subcaptionbox{Constraint between 12:00 and 14:00}{\begin{tikzpicture}[scale=1]
\begin{axis}
[yticklabel=\empty,legend style={at={(0.45,1.23)}, anchor=north,draw=none},legend columns=4,grid=major,xtick  = {0,4,8,12,16,20,24},
width=0.39\textwidth,height=0.3\textwidth,
legend entries={Nominal \quad , Consumption constraint \quad , MPC \space\space\space\space\space\space\space\space\space\space\space\space\space\space\space\space\space\space\space\space\space\space\space\space\space\space\space\space\space\space\space\space\space\space\space\space\space\space\space\space\space\space\space\space\space\space\space\space\space\space\space\space\space\space\space\space\space\space\space\space\space\space\space\space\space\space\space\space\space\space\space}],
\addplot [draw=blue,ultra thick] table[x index=0,y index=1]{DataFigures/nominal.txt};
\addlegendimage{red,ultra thick}
\addplot [draw=Tan,ultra thick, name path=f] table[x index=0,y index=1]{DataFigures/MPC_006_12_14.txt};
\addplot [draw=red,ultra thick, name path=f] table[x index=0,y index=1]{DataFigures/single_constraint_short.txt};
\path [name path=xaxis]
      (\pgfkeysvalueof{/pgfplots/xmin},0) --
      (\pgfkeysvalueof{/pgfplots/xmax},0);
\end{axis}
\end{tikzpicture}}
\hspace{-5.8cm}
\subcaptionbox{Constraint with gradient control}{\begin{tikzpicture}[scale=1]%
\begin{axis}
[grid=major,yticklabel=\empty,xtick  = {0,4,8,12,16,20,24},
width=0.39\textwidth,height=0.3\textwidth],
\addplot [draw=blue,ultra thick] table[x index=0,y index=1]{DataFigures/nominal.txt};
\addplot [draw=Tan,ultra thick, name path=f] table[x index=0,y index=1]{DataFigures/MPC_gradient_control_001.txt};
{DataFiguresOld/AggCons100.0.txt};
\path [name path=xaxis]
      (\pgfkeysvalueof{/pgfplots/xmin},0) --
      (\pgfkeysvalueof{/pgfplots/xmax},0);
\end{axis}
\end{tikzpicture}}
\vspace{-0.15cm}
\caption{(a) Optimized consumption compared to the nominal to limit the global consumption under $30\%$ of the maximum consumption; (b) Limit of $6\%$ between 12:00 and 14:00; (c) Gradient control for the whole day.}
\vspace{-0.2cm}
\label{fig:constraints}
\end{figure*}
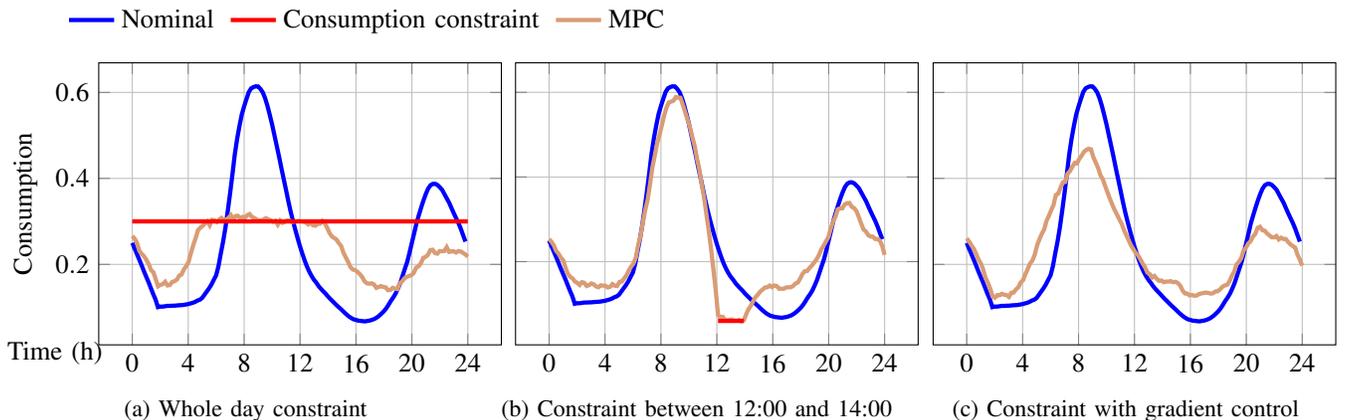

In the first one, we constrain the consumption to be lower than $30\%$ of the maximum consumption, during the whole day. This is done by using the constraint in (\ref{e:upperb}) with $\{u_t\}_{1\leq t \leq T}$ constant $= 0.3$.
This can be especially useful when electric bills incorporate the highest peak in the day. The MPC result is shown in Fig. \ref{fig:constraints}(a).

In the second one (Fig. \ref{fig:constraints}(b)), we try to reduce the consumption drastically during a short window of 2 hours. This can be used to reduce the stress on the electric grid during peak hours.

Finally, to avoid large, sudden changes in consumption, we can apply ramp constraints, as shown in Fig. \ref{fig:constraints}(c). This is done using (\ref{e:gradconstraints}) by setting upper and lower bounds on the gradient of the consumption.

\textbf{In the heterogeneous tracking problem}, we suppose that the WHs differ in volume of water contained (m$^3$), height (m), isolation thickness (m) and heat resistance power (W). As we can see in Fig. \ref{fig:Heterogeneous}, the tracking stays as good.

\subsection{Online Tracking}

We further test the performance of the algorithm to solve a more general problem. In the online setting, the tracking signal is not known in advance. It is only given in real-time, so at time step $t$, we only know the target consumption for the time step $t+\delta t$.

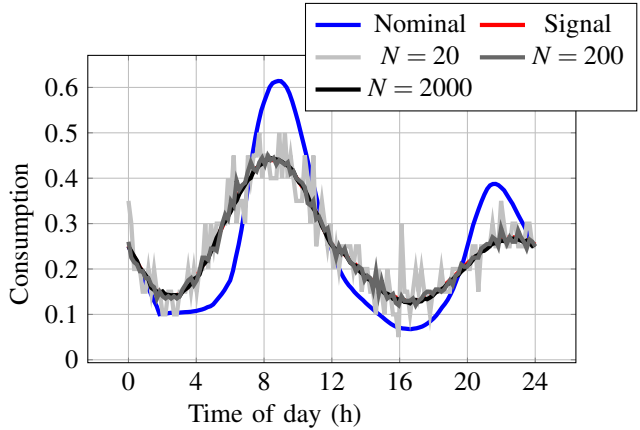
\begin{figure}[ht]
    \centering
    \begin{tikzpicture}[scale=0.45]
    \begin{axis}[xtick  = {0,4,8,12,16,20,24},legend style={at={(2.4,2.5)}, anchor=north east},legend columns=2,
xlabel={Time of day (h)},x label style ={at={(0.75,-0.33)},anchor=north},
grid=major,
width=0.9\textwidth,height=0.6\textwidth,
ylabel={Consumption},y label style ={at={(-0.5,0.25)},anchor=north west},legend entries={Nominal,Signal,$N=20$,$N=200$,$N=2000$}],
    \addplot [draw=blue,ultra thick] table[x index=0,y index=1]{DataFigures/nominal.txt};
    \addplot [draw=red,ultra thick] table[x index=0,y index=1]{DataFigures/Signal_to_track.txt};
    \addplot [draw=white!30!gray!70,ultra thick] table[x index=0,y index=1]{DataFigures/online_tracking_20.txt};
    \addlegendimage{gray!30!black!70,ultra thick}
    \addlegendimage{black,ultra thick}
    \addplot [draw=black,ultra thick] table[x index=0,y index=1]{DataFigures/online_tracking.txt};
    \addplot [draw=gray!30!black!70,ultra thick] table[x index=0,y index=1]{DataFigures/online_tracking_200_better.txt};
    \end{axis}
    \end{tikzpicture}
    \caption{Online tracking setting, with the signal discovered in real-time. 2 switches maximum per WH for $N=2000$ WHs. An average of 3.35 switches per WH for $N=200$ WHs. An average of 8 switches per WH for $N=20$ WHs}
    \label{fig:Online}
\end{figure}

In this setting, the only difference with the algorithm \ref{a:MPC} is that $\lambda$ is a scalar value, and its update in line \ref{a:Evaluation:lbda} of algorithm \ref{a:Evaluation} is as follows:
$$ \lambda_{k}\gets\lambda_{k-1}-\rho_kG_{k,t}.$$

Another such lambda could have been introduced in the MPC algorithm to give more weight to the deviation from the signal at time step $t+\delta t$. This can be introduced to improve the MPC results.

In the online setting, MCOT gives outstanding results, as can be seen in Fig. \ref{fig:Online}. We chose to keep the constraint of a maximum of two switches for $N=2000$ WHs and we relaxed it for $N=200$ and $N=20$.

\section{FUTURE WORK}

The promising results in the heterogeneous case encourage us to explore more complex systems, incorporating various types of flexible loads (such as electric vehicles or air conditioners), while maintaining a global objective.
Another future research direction would be to look at individual costs. In this article, we have looked at a relatively simple cost, but a more in-depth analysis would be interesting, in particular, to see the concrete impact on TCLs of choosing a particular cost function, and to assess the impact on overall constraint compliance performance. Electricity cost could also be incorporated to study different use cases that take into account the global limits of a set of flexible agents (EVs recharging parking lot, etc.), while taking into account evolving electricity price.

\end{document}